\documentclass[journal]{IEEEtran}
\ifCLASSINFOpdf

\else

\fi

\usepackage{amsmath, amssymb}
\usepackage{graphicx}
\usepackage[noadjust]{cite}
\usepackage{color}
\usepackage[caption=false]{subfig}
\usepackage{multirow}
\begin{document}
%An Invariant Set-based Transient Stability-Constrained Optimal Power Flow

\title{A Stability-constrained Optimization Framework for Lur'e Systems with Applications in Power Grids}

\author{Qifeng~Li,~\IEEEmembership{Member,~IEEE,}
         Qiushi~Wang,~\IEEEmembership{Member,~IEEE,} and Konstantin~Turitsyn,~\IEEEmembership{Member,~IEEE}
        % <-this % stops a space
\thanks{Q. Li and K. Turitsyn are with the Department
of Mechanical Engineering, Massachusetts Institute of Technology, Cambridge,
MA, 02139 USA, e-mail: \{qifengli,turitsyn\}@mit.edu;}
\thanks{Q. Wang was with  the Department of Electrical Engineering, Arizona State University, Tempe, AZ 85281. She is now with the American Electric Power, New Albany, OH, 43054, e-mail: qwang@aep.com.}}

\markboth{}%
{Shell \MakeLowercase{\textit{et al.}}: }
\maketitle

\begin{abstract}
For many nonlinear control systems, the chosen equilibrium determines both the steady-state efficiency and the dynamic performance. This paper addresses the issue of obtaining an optimal equilibrium in terms of some steady-state operation criteria for a Lur'e-type system and such an equilibrium can also guarantee a sufficiently large stability region in the dynamic domain such that the system can tolerate some given disturbance. For this purpose, a set of computationally tractable algebraic constraints, which can properly represent the stability certificate under the optimization framework, are proposed. The existing methods formulate the dynamic performance under the optimization framework by discretizing the differential-algebraic equations, which are computationally intractable for large-scale Lur'e systems like power grids. Dissimilarly, the introduced approach first constructs the stability region based on quadratic Lyapunov functions. Then, a novel method is proposed to project the stability region onto the feasible domain of the optimization problem such that the stability certificate can be incorporated into the optimization framework easily. In the transient stability-constrained optimal power flow (TSCOPF) problem of power systems, researchers look for a steady-state operating point with the minimum generation costs that can maintain system stability under some given transient disturbances. The proposed approach is applied to develop a scalable TSCOPF framework for power systems. The TSCOPF model is tested on the IEEE 118-Bus power system.

\end{abstract}


\IEEEpeerreviewmaketitle

\section{Introduction}

\IEEEPARstart{S}{tability} and optimality are two important aspects in the operation and design of a control system. Take the synchronous electric power transmission grid for example: the system operators, on one hand, seek the generation schedule of power plants that meets the demand with the lowest cost; on the other hand, they need to guarantee that the system can maintain synchronism following a typical disturbance. To minimize the costs, the system operators solve an optimal power flow (OPF) problem which is generally a nonlinear programming (NLP) problem \cite{OPF}. At the same time, power system lacks global stability, so the synchrony cannot be a-priori guaranteed for all the possible faults. The outcome of the post-fault dynamics is determined by the solution to a set of differential algebraic equations (DAE) \cite{Kundur}. The optimal dispatch corresponding to the OPF solution determines the system operating point, i.e. its dynamic equilibrium, and it has a significant impact on the system's ability to maintain synchronism after the faults. While currently the OPF and stability assessments are performed separately by System Operators, the next generation of decision support tools should incorporate synchronous stability constraints into the the OPF framework \cite{Scala}.

Unfortunately, the direct incorporation of the dynamic model of power systems into the OPF framework results in a DAE-constrained optimization problem which cannot be solved directly. The most commonly used solution methods for such optimization problems are the discretize-then-optimize approaches \cite{Betts}. Namely, the basic idea of these methods is to discretize the DAE into a set of algebraic equations with respect to small time steps. However, even a small DAE-constrained problem induces a large-scale NLP problem after discretization \cite{Scala,Mak,Hijazi}. As a result, these methods are not practical for large-scale systems like modern transmission grids. Similar challenges are faced by many control systems in other domains, such as Hopfield-type neural networks \cite{VanDen}, machining systems \cite{Lars}, $\Sigma$-$\Delta$ modulators \cite{Abellan}, and flow shop systems \cite{Gokbayrak}. Under the premise of ensuring stability, it is desirable to achieve the best performance of these control systems with respect to reducing costs or increasing efficiency. 

To obtain such a stability-guaranteed optimal solution in a computationally tractable way, this paper proposes a novel framework of stability-constrained optimization for Lur'e-type dynamic systems. The proposed approach is substantially different from the discretize-then-optimize methods and is theoretically scalable. Based on a quadratic Lyapunov function, we cast the stability certificate of Lur'e systems into algebraic constraints in the steady-state domain using LaSalle's invariance principle \cite{Cheng}. The quadratic Lyapunov functions are obtained through the emerging sector nonlinearity approach and the semidefinite programming (SDP) technique \cite{Sector}. Compared with the existing methods, the proposed framework introduces a much smaller computational burden. 

To better illustrate the proposed idea, we introduce an application case in power systems, the transient stability-constrained optimal power flows (TSCOPF) \cite{Gan} and \cite{Scala}. Transient stability is the ability of the power system to maintain synchronism when subjected to a severe transient disturbance \cite{Kundur}. It is heavily affected by the operating point which is the equilibrium of the pre-fault system. TSCOPF is a potentially powerful approach for finding  OPF solutions that can guarantee stability when the system suffers a transient disturbance. Aggressive introduction of renewable generation increases the overall stress of the power system \cite{Li1}, so the stability constraints will likely become the main barrier for transition to clean energy sources. Despite many decades of research, stability assessment is still the most computationally intensive task in power grid operation process. It is even more computationally intractable when stability constraints are considered in optimization problems. This is the motivation for applying the proposed approach to develop a computationally tractable framework for TSCOPF.

The rest of the paper is organized as follows: Section II provides the mathematical formulations of the problem we encounter in power grids; based on the general problem formulations, the proposed stability-constrained optimization framework is presented in Section III; Section IV discusses the application of the proposed approach in TSCOPF of power grids; the novel TSCOPF framework is tested on the IEEE 118-bus system in Section IV; the novelty and limitations of the proposed approach as well as the future research are discussed in detail in the last section.

\section{An Engineering Problem: TSCOPF}

\subsection{Steady-state Model of Power Grids}
The power system is the largest machine in the world composed of a large number of generators and loads interacting through electric flows. Usually, a simplified model of a power transmission grid, which comprises generators, a transmission network, and aggregate loads, is considered for research purposes. Let $\mathcal{N}_G$, $\mathcal{N}_L$, $\mathcal{N}$, and $\mathcal{E}$ denote the sets of generator buses, load buses, all buses, and all edges of the network respectively, the steady-state power network model, i.e. power flow model, is given by
\begin{subequations} \label{PF}
\begin{gather}
\mathbf{\emph{p}}^G-\mathbf{\emph{p}}^L-g^p (\mathbf{\emph{V}},\mathbf{\theta},Y) = 0  \\
  \mathbf{\emph{q}}^G-\mathbf{\emph{q}}^L-g^q (\mathbf{\emph{V}},\mathbf{\theta},Y) = 0 \\
    S(\mathbf{\emph{V}},\mathbf{\theta},Y) \le \overline{S} \\
   \underline{V} \le V \le \overline{V} \\
   \underline{E\theta} \le E\theta \le \overline{E\theta},
\end{gather}
\end{subequations}
where 
\begin{align}
g_i^p(\mathbf{\emph{V}},\mathbf{\theta},Y) &=  V_i\sum_jV_j(G_{ij}\mathrm{cos}(\theta_i-\theta_j)+B_{ij}\mathrm{sin}(\theta_i-\theta_j)) \nonumber \\
&= V_i\sum_jV_j|Y_{ij}|\mathrm{sin}(\theta_{ij}+\alpha_{ij}), \nonumber \\ 
g_i^q(\mathbf{\emph{V}},\mathbf{\theta},Y) &= V_i\sum_jV_j(G_{ij}\mathrm{sin}(\theta_i-\theta_j)-B_{ij}\mathrm{cos}(\theta_i-\theta_j)), \nonumber \\
S_{ij}(\mathbf{\emph{V}},Y)&=Y_{ij}^2V_i^2V_j^2 \quad (i,j \in \mathcal{N},ij \in \mathcal{E}). \nonumber
\end{align}
 
In the above steady-state model, $p_i^G(q_i^G)$ is the vector of active (reactive) generation at bus $i \in \mathcal{N}_G$; $p_i^L(q_i^L)$ is the vector of active (reactive) load at bus $i \in \mathcal{N}_L$; $V$ and $\theta$ are the vectors of bus voltage magnitudes and phase angles respectively in the steady-state domain; $Y$ is a set of alterable parameters and $Y_{ij}$ = $G_{ij}+jB_{ij}$; $G_{ij}$ and $B_{ij}$ are the conductance and susceptance respectively, and $\alpha_{ij} = \arctan(G_{ij}/B_{ij})$; $S_{ij}$ represents the square of apparent power in transmission line $ij$; $E$ is an incidence matrix such that constraint (1e) means $\underline{\theta}_{ij} \le \theta_{ij} \le \overline{\theta}_{ij}$. Note that the rotor velocity $\omega$ of generator is assumed to be constant in steady domain with a uniform value of nearly 1 per unit. Consequently, it is not explicitly included as a variable in the steady-state power flow model.

\subsection{Classical Dynamic Model of Power Grids}

This subsection introduces a dynamic model of power systems which is a standard network-preserving model \cite{Bergen} with transfer conductance included. The classical model is used to formulate the dynamic behaviors of generators in the standard network-preserving model where the reactive power flow is neglected and the magnitude of bus voltage is considered constant during transient. Such a dynamic model may be considered crude for transient stability analysis (TSA). However, it is still computationally unacceptable to incorporate a higher-order dynamic model into the OPF framework.

Let $\delta$ and $\omega$ denote the vectors of bus phase angles and rotor velocities respectively. They are functions of time and $\delta_i(t_0)$=$\theta_i$ ($i \in \mathcal{N}$), $\omega_i(t_0)\approx 1$ per unit ($i \in \mathcal{N}_G$). The dynamic model of angle stability is given as
\begin{subequations} \label{swing}
\begin{align} 
\dot{\delta_i}&=\omega_i-1, \  i \in \mathcal{N}_G \\
   d_i \dot{\delta_i} &= -p^L_i-g_i^p(\mathbf{\emph{V}},\mathbf{\delta},Y), \  i \in \mathcal{N}_L\\
 m_i \dot{\omega_i} &= p_i^G-g_i^p(\mathbf{\emph{V}},\mathbf{\delta},Y)-d_i (\omega_i-1), \  i \in \mathcal{N}_G 
\end{align}
\end{subequations}
where $m_i$ denotes the generator moment of inertia and $d_i$ represents the damping coefficient of generator or load.

The transient disturbances that a power system may encounter include faults on transmission facilities, loss of generation, and loss of large loads \cite{Kundur}. Generally, the disturbance will be cleared after a short period. Hence, from a mathematical perspective, disturbances can be described as the variation of parameters, of which the details are given by
\begin{subequations} \label{fault}
\begin{align} 
(G_{ij}, B_{ij})&=\begin{cases}
          (G_{ij}, B_{ij}) \qquad & t=t_0^-  \\
          (G_{ij}^{\prime\prime}, B_{ij}^{\prime\prime}) \qquad & t=t_0^+ - t_c^-\\
          (G_{ij}^\prime, B_{ij}^\prime) \qquad & t=t_c^+ - \infty
        \end{cases} \\
  p_i^G &= \begin{cases}
    p_i^G \quad (t=t_0^-,t_c^+\rightarrow \infty) \\
    0 \quad (t=t_0^+ \rightarrow t_c^-)
    \end{cases} \\
      p_i^L &= \begin{cases}
    p_i^L \quad (t=t_0^-,t_c^+\rightarrow \infty) \\
    0 \quad (t=t_0^+ \rightarrow t_c^-)
    \end{cases},
\end{align}
\end{subequations}
where $t_c^-=\infty$ if the fault is permanent. An illustrative example of line to ground fault is give in Figure 1. In the pre-fault system, Bus 1 and Bus 2 are connected by a double circuit transmission line and the admittance matrix is $Y$. Suppose, at $t_0$, that a line-to-ground fault is applied at the middle of Line 2, both buses are grounded through half of Line 2. The admittance matrix becomes $Y^{\prime\prime}$. The circuit breakers operate to clear Line 2 at $t_c$ with only Line 1 left. The admittance matrix changes again into $Y^\prime$. Let ($V^\prime,\,\theta^\prime$) denotes the post-fault equilibrium point, the system trajectories are usually required to stay within the polytope $\mathcal{P}=\{\delta_{ij}|-\pi \le \delta_{ij}-\theta_{ij}^\prime \le \pi \}$ during transient.
\begin{figure}[h]
\centering
\includegraphics[width=0.49\textwidth]{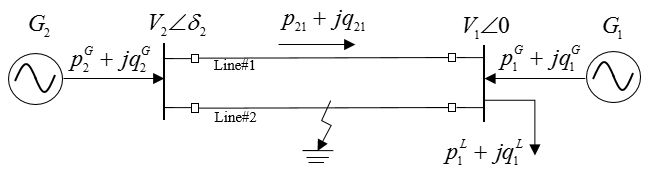}
\caption{The 2-generator system.}
\label{fig:2bus}
\end{figure}

In transient stability assessment (TSA) of power systems, fault type, location, and clearing-time are pre-determined factors. Hence, what determines the system stability is the initial state and the fault-on trajectories are non-trivial functions of the initial point. The fault can generally be cleared within 6 cycles (i.e. 0.1 s) with high-speed relays. Due to the fact that $t_c$ is very small, the non-trivial fault-on trajectories of power grids can be approximated by the following $3^{\text{rd}}$-order Taylor's series with high-fidelity
\begin{subequations}\label{approximation}
\begin{align}
\delta_i(t_c)&=\theta_i-\frac{d_it_c^2}{2!M_i}K_i, i \in \mathcal{N}_G\\
\delta_i(t_c)&=\theta_i-\frac{2!+t_c^2}{2!d_i}K_i, i \in \mathcal{N}_L\\
\omega_i(t_c)&=1+(\frac{d_it_c^2}{2!M_i^2}-\frac{t_c}{M_i})K_i, i \in \mathcal{N}_G,
\end{align}
\end{subequations}
where $K_i = V_i\sum_jV_j(\Delta B_{ij}\mathrm{sin}\theta_{ij}+\Delta G_{ij}\mathrm{cos}\theta_{ij})$, $\Delta G_{ij}=G_{ij}^{\prime\prime}-G_{ij}^\prime$, and $\Delta B_{ij}=B_{ij}^{\prime\prime}-B_{ij}^\prime$. Note that, for the sake of simplicity, the $3^{\text{rd}}$-order term in (4b) is also omitted.

\subsection{Traditional OPF}
A standard formulation of the conventional OPF problem can be expressed as
\begin{align}
\begin{split}
\min_{p^{G}}\quad&f =\sum_{i \in \mathcal{N}_G} (a_{i1}(p_i^G)^2+a_{i2}p_i^G)\\
  \mathrm{s.t.}\quad&\text{(\ref{PF})}.
\end{split} \tag{OPF}
\end{align}
where $f$ is the quadratic cost function of generators. The purpose of (OPF) is to search for a solution of ($p^G$, $q^G$), which can meet the given demand ($p^L$, $q^L$) being subject to the network constraints, with the function $f$ minimized. The solution of (OPF) provides an initial condition to the dynamic system (\ref{swing}) of which the dynamics is significantly affected by this initial condition. However, no information about the dynamic behaviors of (\ref{swing}) is considered in (OPF). As a result, the initial condition obtained by solving (OPF) can not guarantee the stability when a given fault occurs.

Based on the 2-generator system given in Figure \ref{fig:2bus}, we provide an example
to illustrate the stability issue with model (OPF). Suppose that the system is lossless (i.e. the resistance in the transmission line is ignored) and the generation costs $c_1 > c_2$. The result of solving (OPF) model for this system is $p^G_2=p_{21}=p^L_1$ and $p^G_1=0$. Further assume that Bus 1 is an infinite bus which has sufficient capability to keep $V_1$ constant. The result of transient stability assessment shows that, to maintain the system's stability under the line fault on Line$\#$2, there is a certain upper limit on the line power $p_{21}$, i.e. $p_{21} < \overline{p}_{21} < p^L_1$. A stability-guaranteed optimal solution should be $p^G_2=\overline{p}_{21}$ and $p^G_1=p^L_1-\overline{p}_{21}$. However, the model (OPF) fails to include the quantity $\overline{p}_{21}$ which is directly related to transient stability of the 2-bus system.

\subsection{State-of-the-Art of TSCOPF}
TSCOPF has been receiving a growing amount of attention as a potential solution to the problem introduced in the above subsection. The following TSCOPF model is widely adopted in literatures, such as \cite{Gan}, \cite{Scala}, \cite{Chen}, and \cite{Conejo}, 
\begin{align}
\begin{split}
\min_{p^{G}}\quad&f =\sum_{i \in \mathcal{N}_G} (a_{i1}(p_i^G)^2+a_{i2}p_i^G)\\
  \mathrm{s.t.}\quad&\text{(\ref{PF}), (\ref{swing}), (\ref{fault})}, and \\
  &\delta(t) \in \mathcal{P} \quad \forall t.
\end{split} \tag{TSCOPF1}
\end{align}
However, the model (TSCOPF1) has two evident drawbacks: i) transient-stability can not be strictly guaranteed by just restricting the system trajectories $\delta(t)$ in a given polytope; ii) (TSCOPF1) can not be solve directly due to the DAE constrains (\ref{swing}). To solve such a DAE-constrained optimization problem, former researchers discretized the DAEs in (\ref{swing}) into a set of algebraic constrains in terms of small time steps. However, it results in a large-scale nonlinear programming (NLP) problem which is computationally intractable for large-scale system like power grids. 

Since (TSCOPF1) is not practical, power grid operators currently use iterative OPF algorithms where an independent TSA is required at each iteration \cite{Cai} and \cite{Pizano}. Note that the dynamic generator model used in literature is generally the classical model. To further reduce the computational burden, some existing methods are based on simpler dynamic models of power systems like the single machine equivalent model used in \cite{Xia}. However, these approaches are still not practical in real-time application due to the heavy computational burden they incur. Our former work shows that it is possible to reformulate the dynamic power system model, which is described by the classical generator model with transfer conductance considered \cite{Long1} and network preserved \cite{Long2}, into the classical Lur'e form. Next section proposed a general stability-constrained optimization framework for Lur’e systems. Based on this novel optimization framework, we will develop a computationally tractable model for TSCOPF.

\section{Stability-constrained Optimization Framework for Lur'e Systems}
Motivated by the above problem engineers encountered in power systems, this section aims at developing a novel stability-constrained optimization framework which: 1) is applicable to various Lur'e-type control systems, 2) can guarantee strict system stability, and 3) is computationally tractable. This section starts introducing the novel stability-constrained optimization framework with general problem formulations, so that it is convenient to apply the proposed approach to other control systems other than power grids. 

\subsection{A General formulation of the Problem}

\begin{figure}
\centering
\includegraphics[width=0.32\textwidth]{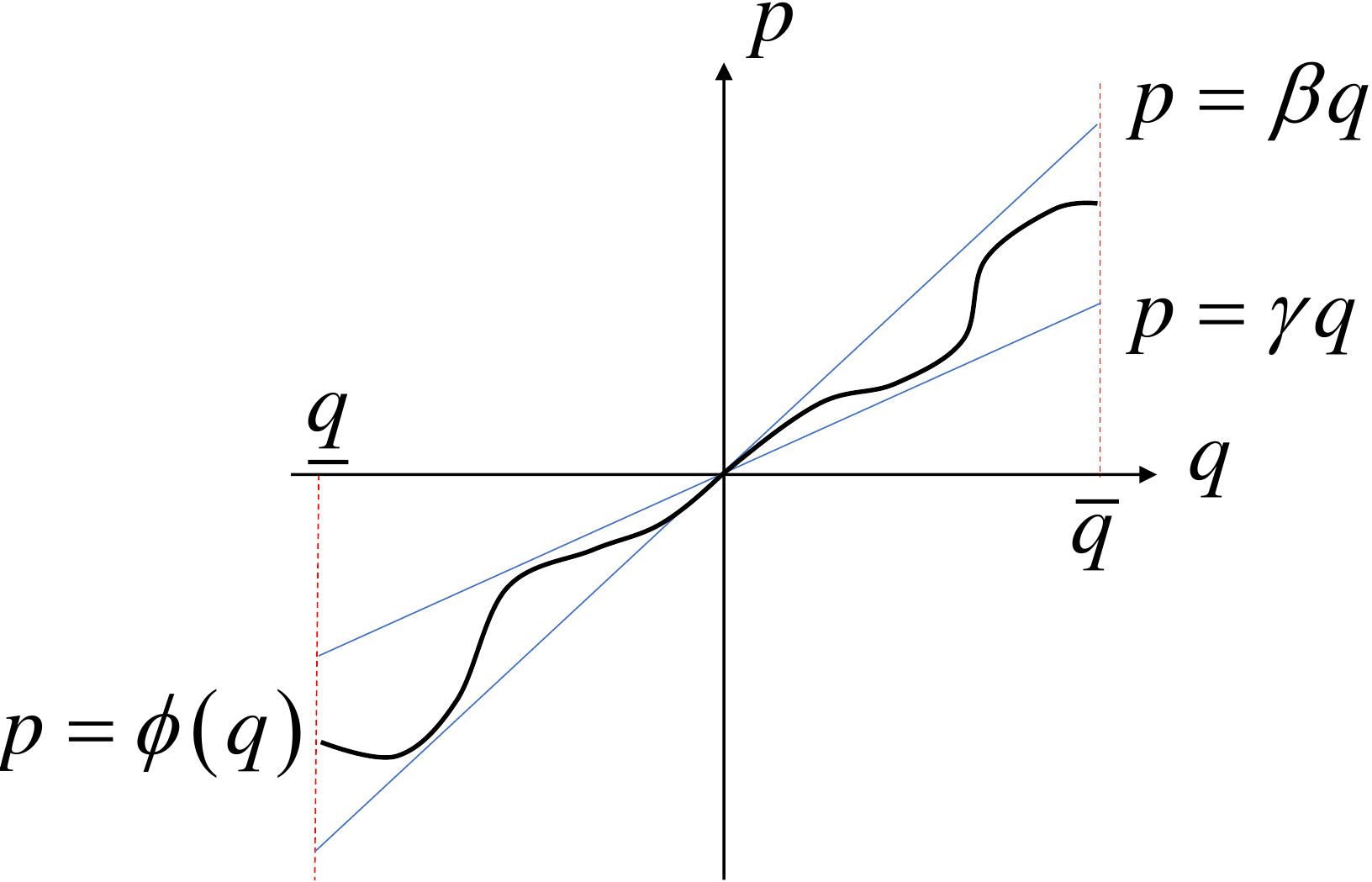}
\caption{The concept of sector bounds.}
\label{fig:sectorbounds}
\end{figure}

Let's consider the following Lur'e-type autonomous system with a nonlinear term of vanishing disturbance
\begin{subequations} \label{Lure}
\begin{align}
 \dot x &= A x + B \phi(y)+ B_u\begin{cases}
    0 \quad (t=t_0^-,t_c^+\rightarrow \infty) \\
    u(x) \quad (t=t_0^+ \rightarrow t_c^-)
    \end{cases}\\
    y&=Cx,
\end{align}
\end{subequations}
where $x$ is a vector function of time $t$ representing the dynamic state vector; the nonlinearity $\phi:\mathcal{R}^n \rightarrow \mathcal{R}^m$ can be bounded by some local sectors as shown in Fig \ref{fig:sectorbounds} and the sector is defined below. In the above formulation, we define the dynamic state vector $x$ in the coordinates where the equilibrium point of the system is the origin, which means $x(t_0)=0$ and $\phi(0)=0$. System (\ref{Lure}) can also be considered as a linear system with nonlinear feedback $\phi$ and an exogenous input $u$. Take the case in power grids for example. We can observe by comparing (1a) and (\ref{swing}) that the angle dynamics model can be reformulated into the form of (\ref{Lure}) if defining $x=[\delta-\theta^\prime;\, \omega-\omega_0]$. Please refer to Subsection IV-A for more details about this reformulation.

\textbf{Definition.} \textit{The nonlinear function $\phi$ is said to be locally bounded by sector} [$\gamma,\beta$] \textit{if, for all $q \in$} [$\underline{q},\overline{q}$]\textit{, $p = \phi(q)$ lies between $\gamma$ and $\beta$}. 

In steady-state, operators generally have freedom to adjust the control vector $z$ to achieve the best system performance in terms of the criterion $f$ by solving the following optimization problem
\begin{subequations} \label{Optimization}
\begin{align}
\min_{z} \; &f(x^*,z) \label{Objective}  \\ 
\mathrm{s.t.}\; &g(x^*,z) = 0 \label{NonLinearSyst} \\
&h(x^*,z) \le 0, \label{PhysicalConstraint}
\end{align}
\end{subequations}
where $x^*$ denotes the vector of steady-state variables and the equilibrium of the dynamic system (\ref{Lure}). Between (\ref{Lure}) and (\ref{Optimization}), there exist the following relations
\begin{subequations}
\begin{gather} \label{Connection}
x(t)=\tilde{x}(t)-x^* \  \text{and} \  \tilde{x}(t_0)=x^* \nonumber \\
g(x+x^*,z) =A x - B \phi(C x) \\
u(x)=v(x+x^*)
\end{gather}
\end{subequations}
where $\tilde{x}(t)$ is the vector of dynamic state variables in $x^*$-coordinates and $v(\cdot)$ is a derivable function. Note that the control variable $z$ of the steady-state optimization problem (\ref{Optimization}) is a fix parameter in dynamics and, therefore, not explicitly included in the dynamic formulation (\ref{Lure}). In a power transmission system, constraints (\ref{NonLinearSyst}) and (\ref{PhysicalConstraint}) correspond to (1a)-(1b) and (1c)-(1e) respectively if we let $x^*=[V,\,\theta]^T$ and $z = [p^G,q^G]^T$.

\subsection{Scalable Algebraic Stability Certificate}

This subsection develops a scalable algebraic certificate of stability which can be incorporated into the optimization model (\ref{Optimization}). Without loss of generality, we assume that the trajectories of system (\ref{Lure}) are required to stay within the polytope $\mathcal{P}=\{x|\underline{y} \le Cx \le \overline{y}\}$ similar to the power system. Let $W(x)= x^TPx$ be a positive definite function. According to the classical circle criterion \cite{Khalil} and the related approaches presented in \cite{Sector}, if the sector [$\gamma,\beta$] is valid for locally bounding the nonlinearity $\phi$ of (\ref{Lure}) in the the polytope $\mathcal{P}$, the positive definite matrix $P$ can be obtained by solving the following linear matrix inequality (LMI) \cite{Boyd}
\begin{equation} \label{LMI}
\left[
\begin{array}{cc} A^TP+PA -   C^T\tau\gamma\beta C     & PB+\frac{1}{2}(\gamma+\beta)C^T \\
 B^TP+\frac{1}{2}(\gamma+\beta)C  & -\tau \\
\end{array}
\right] \preceq 0.
\end{equation}
where $\tau \geq 0$. Moreover, we have the following proposition.

\textbf{Proposition}. \textit{$\dot{W}(x)$ is strictly negative if}
\begin{align} \label{lemma}
(\phi (Cx) - \gamma Cx)^T(\phi (Cx) - \beta Cx) < 0 
\end{align}
\textit{holds for all $x$ in $\mathcal{P}$-$\{0\}$.}

\textit{Proof}: See Appendix A. $\square$

The proposition implies that, by carefully selecting the sector, we are able to obtain a $P$ that makes $W(x)$ strictly decay along the system trajectories in polytope $\mathcal{P}$. Following from LaSalle's \textit{invariance principle}, we have the lemma below.

\textbf{Lemma}. \textit{For any fault-cleared state $x(t_c)$ within the set $\Omega$ which is defined by}
\begin{equation} \label{energy}
 \Omega = \{x(t) \in \mathcal{P} \mid W(x(t_c)) \le W^{\min}\},
  \end{equation}
\textit{where} \begin{align} \label{Wmindefinition}
W^{\min}=\min_{x\neq 0}\{W(x)|\ Cx=\underline{y},\ \text{or}\  Cx=\overline{y}\}, 
\end{align}
\textit{the system trajectories of the nominal system (\ref{Lure}) starting from $x(t_c)$ stay in the set $\Omega$ for all $t \geq t_c$ and eventually will converge to the origin.}

\textit{Proof}: See Appendix B. $\square$

Note that, any feasible point $x^*$ of (\ref{Optimization}) is a potential equilibrium of the nominal system of (\ref{Lure}). Since the exact equilibrium is unknown before solving the optimization problem, it is necessary to choose a sector that is valid for all the potential equilibria, namely the whole steady-state feasible set of (\ref{Optimization}), rather than just one single equilibrium. To obtain such a sector, one needs to regard $x^*$ as a changeable parameter in $\phi$ and consider its range of variation. An example of selecting such a qualified sector for the case in power systems is given in Figure \ref{fig:sectors}. Moreover, if $\phi(q)$ = [$\phi_1(q_1),\dots,\phi_k(q_k)$]$^T$ represents multiple nonlinearities, one can customize the chose sector for each nonlinearity such that $\beta=\emph{\emph{diag}}(\beta_1,\dots,\beta_k)$, $\gamma=\emph{\emph{diag}}(\gamma_1,\dots,\gamma_k)$, and $\tau=\emph{\emph{diag}}(\tau_1,\dots,\tau_k)$.

Even though $W(x)$ is quadratic, it is still computationally intractable to search for $W^{\min}$ of large-scale systems through (\ref{Wmindefinition}) under the optimization framework. We can equivalently rewrite (\ref{Wmindefinition}) as
\begin{equation}
\mathit{W}^{\min}=\min_{i}\mathit{W}_i^{\min}  \nonumber
\end{equation}
\begin{equation} \label{eq4}
\begin{split}
 \mathit{W}_i^{\min}=\min_{x}x^TPx \\
\mathrm{s.t.}\quad C_i^Tx=y_i
\end{split}
\end{equation}
where $C_i^T$ is the $i$-th row of matrix $C$, and $y_i$ = $\overline{y}_i$ or $\underline{y}_i$. According to the first-order optimality conditions \cite{NumericalOptimization}, the optimization problem (\ref{eq4}) has a trivial solution $\hat{x} = y_i P^{-1} C_i/(C_i^T P^{-1}
C_i)$. Consequently, $\mathit{W}^{\min}$ can be obtained through the following simpler way
\begin{equation}
\mathit{W}^{\min}=\min_{i}\{\frac{\min\left\{\overline{y}_i^2,\underline{y}_i^2\right\}}{C_i^T
P^{-1} C_i}\}. \nonumber
\end{equation}
Generally, in the $x$-coordinates, $\overline{y}_i=\Delta l$ and $\underline{y}_i=-\Delta l$, where $\Delta l$ is positive, while  $\overline{y}_i=C_i^Tx^* +\Delta l$ and $\underline{y}_i=C_i^Tx^* -\Delta l$ in the $x^*$-coordinates. Hence, we can search for $\mathit{W}^{\min}$ in the $x^*$-coordinates though 
\begin{equation}  \label{Wmin}
\mathit{W}^{\min}=\min_{i}\{\frac{\min\left\{(C_i^Tx^* +\Delta l)^2,(C_i^Tx^* -\Delta l)^2\right\}}{C_i^T
P^{-1} C_i}\}. 
\end{equation}

\subsection{A Scalable Stability-constrained Optimization Framework}

At $t=t_0$, system (\ref{Lure}) suffers a disturbance and its trajectory $x(t)$ starts deviating from the origin (i.e. the equilibrium point in the $x^*$-coordinates). Let $x(t_c)$ denote the state when the disturbance is cleared. As discussed before that the fault clearing time $t_c$ is  short, the fault-cleared state can be accurately approximated via the following Taylor's series 
\begin{align} \label{trajectory}
  x(t_c) &= x(t_0) +\sum_{n=1}^{N} \frac{x^{(n)}(t_0)}{n!}(t_c-t_0)^n
  \nonumber \\
  &=\sum_{n=1}^{N}\frac{g^{(n-1)}(x^*)+v^{(n)}(x^*)}{n!}(t_c-t_0)^n,
\end{align}
where $A x(t_0) - B \phi(C x(t_0))=g(x^*,z)=0$ and $u(x(t_0))=v(x^*)$ as defined in Subsection III-A. As discussed before, $N=3$ is sufficient for power grid transient stability.

According to the Lemma, the post-disturbance trajectory $x(t \ge t_c)$ can stay within $\mathcal{P}$ or even converge back to the origin as $t \rightarrow \infty$ if the fault-clearing state $x(t_c)$ satisfies condition (\ref{energy}). Note that $W^{\min}$ in (\ref{Wmin}) and $x(t_c)$ in (\ref{trajectory}) are functions of $x^*$. Hence, constraints (\ref{energy}) and (\ref{Wmin})-(\ref{trajectory}) together represent the dynamic stability certificate in the $x^*$-domain. By adding (\ref{energy}) and (\ref{Wmin})-(\ref{trajectory}) to problem (\ref{Optimization}), we have the following stability-constrained optimization model for the nonlinear system (\ref{NonLinearSyst})-(\ref{PhysicalConstraint})/(\ref{Lure})
\begin{align}  
\begin{split}
\min_{x} \quad &\text{(\ref{Objective})} \\
\text{s.t.} \quad &\text{(\ref{NonLinearSyst}),\,(\ref{PhysicalConstraint})},\,(\ref{energy}),\,(\ref{Wmin}),\,\text{and}\,(\ref{trajectory})
\end{split} \tag{BL-SCO}.
\end{align}

Model (BL-SCO) is a bilevel optimization problem, where (\ref{Wmin}) is the lower-level subproblem. The stability of (\ref{Lure}) has been taken into account in (BL-SCO). However, it is still very hard to solve since the subproblem (\ref{Wmin}) is nonconvex. To overcome this issue, we proposed the following single-level optimization model 
\begin{align}  
\text{(SL-SCO)}\quad\min_{z} &F(x^*,z,W^{\min}) = f(x^*,z) - \epsilon W^{\min} \label{Objective1} \\
\mathrm{s.t.}\quad &\text{(\ref{NonLinearSyst}),\,(\ref{PhysicalConstraint})}, (\ref{energy}), (\ref{trajectory}), \; \textmd{and} \nonumber \\
&\begin{cases}  \label{Concave}
W^{\min} \le \frac{(C_i^Tx^* - \Delta l)^2}{C_i^T
P^{-1} C_i} \\ 
W^{\min} \le \frac{(C_i^Tx^* + \Delta l)^2}{C_i^T
P^{-1} C_i} 
\end{cases}.
\end{align}

\textbf{Theorem 1}. \textit{Optimal solution of (SL-SCO) is also optimal to optimization problem (BL-SCO).}

\textit{Proof}: See Appendix C. $\square$

\textbf{Remark 1}. Model (SL-SCO) is the novel stability-constrained optimization framework for nonlinear system (\ref{Lure}) proposed in this paper. Theorem 1 implies that (SL-SCO) is equivalent to (BL-SCO). Although the feasible set of (SL-SCO) is just a relaxation of (BL-SCO)'s, by adding the perturbation term $- \epsilon W^{\min}$ to the objective function, one can obtain an exact locally optimal solution of problem (BL-SCO) by solving (SL-SCO).

Problem (SL-SCO) is a single-level optimization problem and much easier to solve than (BL-SCO). Note that both (BL-SCO) and (SL-SCO) are nonconvex. Therefore, only locally optimal solutions can be guaranteed. The feasible region specified by constraints (\ref{energy}), (\ref{trajectory}) and (\ref{Concave}) together can be regarded as the projection of invariant set $\Omega$ in the $x^*$-domain. Inequalities (\ref{Concave}) are concave which belong to a special type of non-convex constraints. A discussion on convexifying (\ref{Concave}) is given in next subsection. 

\subsection{Two Convex Options of $W^{\min}$ Calculation}

\begin{figure}[tb]
\centering
\includegraphics[width=0.48\textwidth]{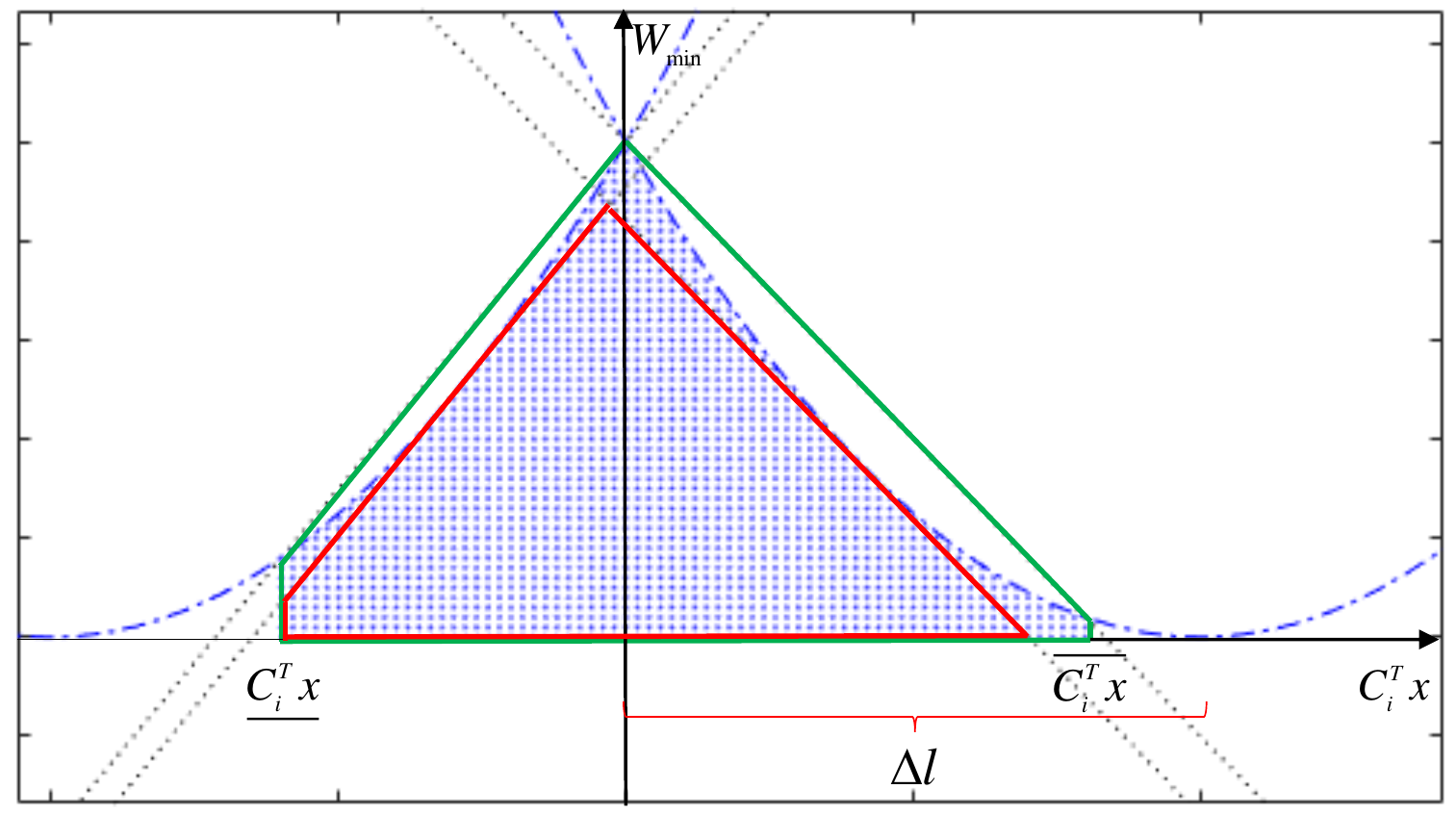}
  \caption{Two convex alternatives of concave constraints (\ref{Concave}). The shaded area denotes the original concave feasible set of (\ref{Concave}). The regions with green and red boundaries are the convex hull relaxation and a convex inner approximation of the shaded area respectively.}
  \label{fig:convexhull}
\end{figure}

Convex optimization, due to its high-efficiency, has been applied to a wide range of automatic control systems \cite{Boyd2}. For example, the convexification of OPF, including convex relaxations \cite{Low} and convex inner approximations \cite{Hung}, is recently one of the research hotspots in the optimization and power sectors. Among the stability constraints in (SL-SCO), constraint (\ref{trajectory}) is related to the system (\ref{NonLinearSyst}) and (\ref{PhysicalConstraint}), while (\ref{energy}) and (\ref{Concave}) are convex and concave respectively. To meet the future need of obtaining a convex stability-constrained optimization framework for the nonlinear system (\ref{NonLinearSyst}) and (\ref{PhysicalConstraint}), this subsection offers two convex alternatives of the concave constraint (\ref{Concave}): a convex hull relaxation and a convex inner approximation. A pictorial interpretation is given in Figure \ref{fig:convexhull}.

\textbf{Theorem 2}. \textit{Set $\Psi$ is the convex hull of set $\psi$, where}
 \[
    \psi = \left\{(x^*,W^{\min}) \left| \begin{array}{lr}
W^{\min} \le \frac{((C_i^Tx^* - \Delta l)^2 }{C_i^T
P^{-1} C_i} \\ 
W^{\min} \le \frac{(C_i^Tx^* + \Delta l)^2}{C_i^T
P^{-1} C_i} \\
\underline{C_i^Tx^*} \le C_i^Tx^* \le \overline{C_i^Tx^*}
\end{array}\right. \right\} 
  \]
   \[   
    \Psi = \left\{(x^*,W^{\min}) \left| \begin{array}{lr}
W^{\min} \le \frac{(\overline{C_i^Tx^*} - 2\Delta l)C_i^Tx^*+\Delta l^2 }{C_i^T
P^{-1} C_i} \\ 
W^{\min} \le \frac{(\underline{C_i^Tx^*} + 2\Delta l)C_i^Tx^*+\Delta l^2}{C_i^T
P^{-1} C_i} \\
\underline{C_i^Tx^*} \le C_i^Tx^* \le \overline{C_i^Tx^*}
\end{array}\right. \right\}, 
  \]
 \textit{where $W^{\min}$ is nonnegative and} $-2\Delta l \le \underline{C_i^Tx^*} \le 0 \le \overline{C_i^Tx^*} \le 2\Delta l $.\\
\textit{Proof}: See Appendix D. $\square$

Replacing the concave constraint (\ref{Concave}) with its convex hull relaxation, we have a relaxed version of the proposed stability-constrained optimization frame (SL-SCO), which is given by
\begin{equation} \label{framework}
\begin{split}
 \min \  &(\ref{Objective1}) \\
\mathrm{s.t.}\  &\text{(\ref{NonLinearSyst}),\,(\ref{PhysicalConstraint})}, (\ref{energy}), (\ref{trajectory}), \; \textmd{and} \, \Psi.
\end{split} \tag{R-SLSCO}
\end{equation}

According to Theorem 2 and the definition of convex hull \cite{Li2}, set $\Psi$ is the best convex relaxation of $\psi$. Compared with the concave constraints in $\psi$, the linear constraints in $\Psi$ are easier to compute. It is a consensus that the Lyapunov-based methods for stability assessment are more or less conservative. Searching for the optimal solutions over set $\Psi$ will result in slightly bigger values for $W^{\min}$, which will reduced the conservativeness of the Lyapunov method. Nevertheless, an additional stability assessment may be required due to the optimism introduced by the relaxation. From an empirical perspective, the convex relaxation version (\ref{framework}) is recommended for power transmission grids.

On the under hand, if (\ref{Concave}) is replaced by the proposed inner approximation, we have a more strictly stable version of (SL-SCO):
\begin{equation} 
\begin{split}
 \min \  &(\ref{Objective1}) \\
\mathrm{s.t.}\  &\text{(\ref{NonLinearSyst}),\,(\ref{PhysicalConstraint})}, (\ref{energy}), (\ref{trajectory}), \; \textmd{and}\\ 
&\begin{cases}  
W^{\min} \le \frac{(\overline{C_i^Tx^*} - 2\Delta l)C_i^Tx^*+\Delta l^2-\overline{C_i^Tx^*}^2/4}{C_i^T
P^{-1} C_i} \\ 
W^{\min} \le \frac{(\underline{C_i^Tx^*} + 2\Delta l)C_i^Tx+\Delta l^2-\underline{C_i^Tx^*}^2/4}{C_i^T
P^{-1} C_i} 
\end{cases}
\end{split} \tag{I-SLSCO}
\end{equation}
where the inner approximation is developed by considering the mean value point on both sides of the origin, i.e. ($\underline{C_i^Tx^*}/2$, $(\underline{C_i^Tx^*}/2+\Delta l)^2$) and ($\overline{C_i^Tx^*}/2$, $(\overline{C_i^Tx^*}/2-\Delta l)^2$), as the tangent points. Note that the above formulation is not the unique option of inner approximations. Actually, the best convex inner approximation of a nonconvex set is yet to be defined. Defining such an inner approximation will be one of the tasks in our future research. For the cases that quadratic Lyapunov functions are not very conservative, the inner approximation version (I-SLSCO) is recommended to strictly guarantee stability.

\section{Application in Power Grids: TSCOPF}

This section applies the framework (\ref{framework}) to the power system introduced in Section II to develop a novel optimization model of TSCOPF. Before introducing the new TSCOPF framework, more details about Lur'e-type reformulation of the angle dynamics in power transmission systems and selection of the non-uniform sectors.
\subsection{Lur'e-type Reformulation of Power Grid Dynamics}

The angle dynamic model of power grids adopted in this paper is an improved network-preserving model where the transfer conductance is retained. An important assumption of the classical generator model is that the voltage $V$ keeps constant during the first swing of transient. In this paper, we use equations (1a) and (1b) to calculate the pre-fault equilibrium ($V,\,\theta$) and post-fault equilibrium ($V^\prime,\,\theta^\prime$) based on the pre-fault conditions ($p^G,\,q^G,\,p^L,\,q^L,\,Y$) and post-fault conditions ($p^{G\prime},\,q^{G\prime},\,p^{L\prime},\,q^{L\prime},\,Y^\prime$) respectively. If the system cannot restore exactly to the pre-fault state after the fault is cleared, we have $V^\prime \neq V$ and $\theta^\prime\neq \theta$. This section introduces the details of reformulating the angle dynamic model (\ref{swing}) - (\ref{fault}) of transmission grids into the form of (\ref{Lure}).

Considering the post-fault system as the nominal system and the post-fault equilibrium point as the origin, the dynamic-state variable $x$ is defined as $x = [\delta_1 -
\theta_1^\prime,\dots, \delta_{|\mathcal{N}_G|} -
\theta_{|\mathcal{N}_G|}^\prime, \delta_{|\mathcal{N}_G|+1} -
\theta_{|\mathcal{N}_G|+1}^\prime,\dots, \delta_{|\mathcal{N}|} -
\theta_{|\mathcal{N}|}^\prime,\omega_1-1,\dots,\omega_{|\mathcal{N}_G|}-1]^T$. Note that, in steady-state (including the pre-fault and post-fault steady state), the system has a uniform angle velocity which is nearly 120$\pi$ rad/s or 1 in per unit. Let $E$ be the incidence matrix of the directed graph $\mathcal{G}(\mathcal{N},\mathcal{E})$, so
that $E[\delta_1,\dots,\delta_{|\mathcal{N}|}]^T =
[(\delta_i-\delta_j)_{\{i,j\}\in\mathcal{E}}]^T$. Let matrix $C=E[I_{|\mathcal{N}| \times |\mathcal{N}|} \;O_{|\mathcal{N}| \times |\mathcal{N}_G|}]$, then 
$$Cx=E[\delta_1 -
\theta_1^\prime,\dots, \delta_{|\mathcal{N}|} -
\theta_{|\mathcal{N}|}^\prime]^T=[(\delta_{ij}-\theta_{ij}^\prime)_{\{i,j\}\in\mathcal{E}}]^T.$$

For dynamic model (\ref{swing}), the matrices $A$ and $B$ in (\ref{Lure}) are given as 
$$A=\left[
        \begin{array}{ccccc}
          O_{|\mathcal{N}_G| \times |\mathcal{N}|} \qquad & I_{|\mathcal{N}_G| \times |\mathcal{N}_G|}\\
          O_{|\mathcal{N}|-|\mathcal{N}_G| \times |\mathcal{N}|} \qquad & O_{|\mathcal{N}|-|\mathcal{N}_G| \times |\mathcal{N}_G|} \\
          O_{|\mathcal{N}_G| \times |\mathcal{N}|} \qquad &-M_1^{-1}D_1
        \end{array}
      \right]$$
and $$
 B= \left[
        \begin{array}{ccccc}
          O_{|\mathcal{N}_G| \times |\mathcal{E}|}; \quad
          S_1M^{-1}E^T; \quad
          S_2M^{-1}E^T
        \end{array}
      \right]$$
respectively, where $S_1$ = $[O_{|\mathcal{N}|-|\mathcal{N}_G| \times |\mathcal{N}_G|}\; I_{|\mathcal{N}|-|\mathcal{N}_G|\times |\mathcal{N}|-|\mathcal{N}_G|}],$ and $S_2$ = $[I_{|\mathcal{N}_G| \times |\mathcal{N}_G|} \, O_{|\mathcal{N}_G|\times |\mathcal{N}|-|\mathcal{N}_G|}]$; $M_1$ = $\emph{\emph{diag}}(m_1,\dots,m_{|\mathcal{N}_G|})$, $D_1=\emph{\emph{diag}}(d_1,\dots,d_{|\mathcal{N}_G|})$ and $M=\emph{\emph{diag}}(d_{|\mathcal{N}_G|+1},\dots,d_{|\mathcal{N}|},m_1,\dots,m_{|\mathcal{N}_G|})$ are the matrices of moment of inertia, frequency controller action on governor, and  frequency coefficient of load respectively.

Following the definitions of matrices $A$ and $B$, we consider the vector of nonlinear interactions $\phi$ in the simple trigonometric form
\begin{equation}
\phi_k(C_k^Tx)=V_i^\prime V_j^\prime|Y_{ij}^\prime|(\sin(\delta_{ij}+\alpha_{ij})-\sin(\theta_{ij}^\prime+\alpha_{ij})), \nonumber
\end{equation}
where $C_k^T$ is the $k$th row of $C$ and the term $V_i^\prime V_j^\prime|Y_{ij}^\prime|$ can be regarded as a parameter for any given post-fault equilibrium during fault-on transient. In next subsection, we will show that, by defining the nonlinear interactions $\phi$ in this way, one can customize sector bounds for each nonlinearity and obtain valid quadratic Lyapunov functions for the whole feasible region rather than a single equilibrium point.

Finally, this paragraph introduces the perturbation term of the angle dynamic model (\ref{swing}) based on the time-variant parameters given in (\ref{fault}). As mentioned above, the stability certificate is constructed based on the post-fault nominal system and, consequently, not directly related to the disturbance. The information of disturbance is only required in the acquisition of fault-cleared state. Therefore, for the cases where the pre- and post-fault conditions are not identical, it is more convenient to observe the disturbance term from the pre-fault equilibrium in time period [$t_0,\,t_c$]. The detailed expression is given as
\begin{align}
&\begin{cases}
          0 \; & t=t_0^-  \\
          V_i\sum_jV_j(\Delta B_{ij}\mathrm{sin}\theta_{ij}+\Delta G_{ij}\mathrm{cos}\theta_{ij}) \; & t=t_0^+ - t_c^-
 \end{cases} \nonumber \\
        &\begin{cases}
    0 \quad &t=t_0^- \\
    -p_i^G \quad &t=t_0^+ \rightarrow t_c^- 
    \end{cases} \nonumber \\
            &\begin{cases}
    0 \quad &t=t_0^-,t_c^+\rightarrow \infty \\
    p_i^L \quad &t=t_0^+ \rightarrow t_c^- 
    \end{cases}, \nonumber
\end{align}
 where $\Delta G_{ij}=G_{ij}^{\prime\prime}-G_{ij}$ and $\Delta B_{ij}=B_{ij}^{\prime\prime}-B_{ij}$.

\subsection{Lyapunov Functions and Fault-on Trajectories}

In traditional transient-stability analysis of power systems, the Lyapunov function is constructed based on a given post-fault equilibrium point \cite{Long1} and \cite{Long2}. However, the Lyapunov function obtained in this way is not valid for the TSCOPF framework since the post-fault equilibrium is a solution of TSCOPF which is unknown before solving the TSCOPF. This subsection aims at constructing an effective common Lyapunov function for the whole feasible region in the steady-state domain by carefully designing the sector bounds.

\begin{figure}[tb!]
\centering
\includegraphics[width=0.45\textwidth]{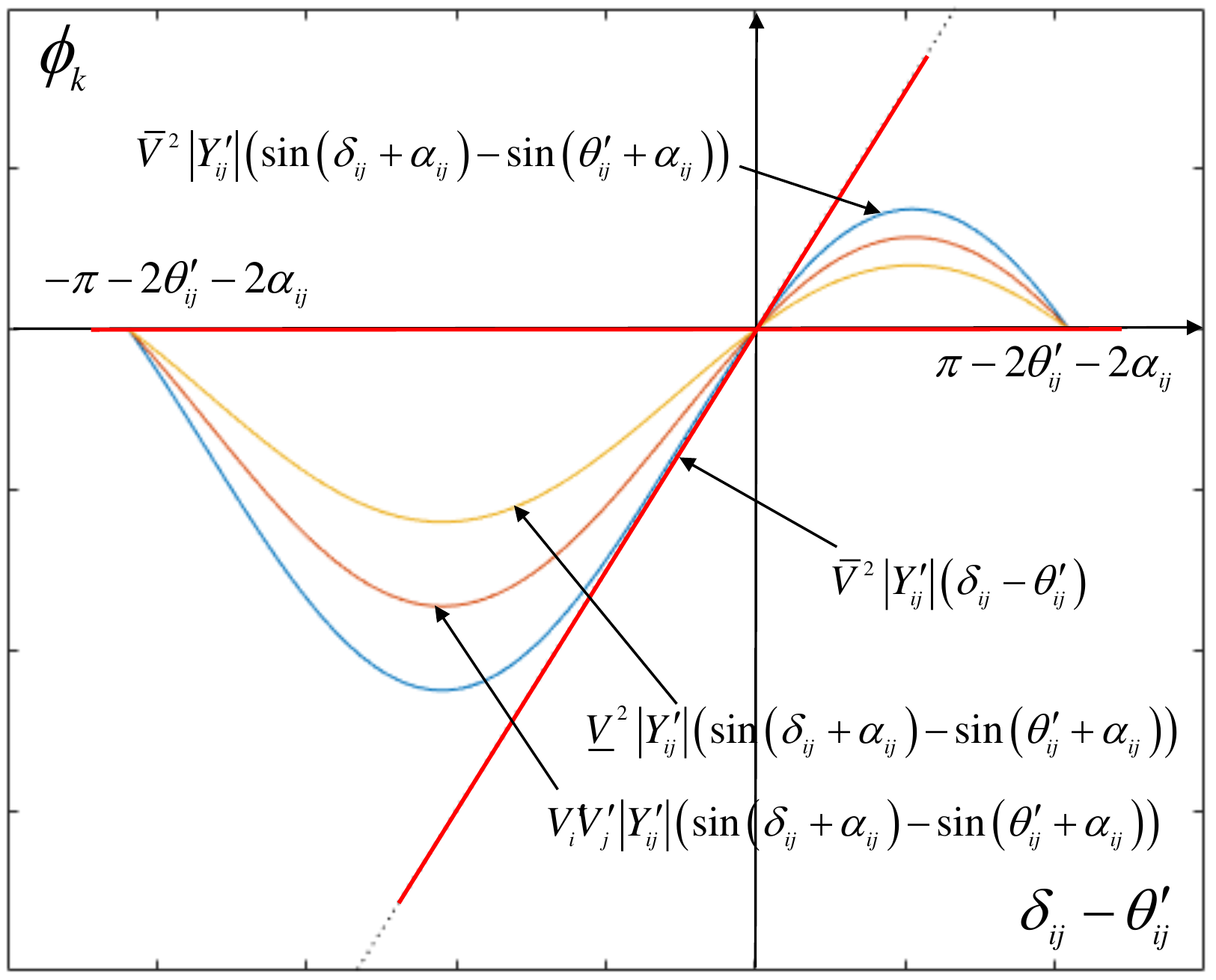}
  \caption{The designed sector bounds for the nonlinearities in system (\ref{swing}).}
  \label{fig:sectors}
\end{figure}

Assuming that $\overline{y}=\pi-2\theta_{ij}^\prime-2\alpha_{ij}$ and $\underline{y}=-\pi-2\theta_{ij}^\prime-2\alpha_{ij}$, we can observe from Fig. \ref{fig:sectors} that, within the polytope $\mathcal{P}$ = $\{\delta_{ij}|-\pi-\theta_{ij}^\prime-2\alpha_{ij}\le \delta_{ij} \le \pi-\theta_{ij}^\prime-2\alpha_{ij}, (i,j)\in\mathcal{E}\}$, the nonlinear terms $\phi_k(C_k^Tx)$ ($k \in \mathcal{E}$ is the serial number corresponding to branch $(i,j)$ in the branch set $\mathcal{E}$) are bounded by the sectors [0, $\beta_k$] for any post-fault equilibrium if $\beta_k=\overline{V}^2|Y_{ij}^{\prime}|$. It has been pointed out in \cite{Long2} that, for the structure-preserving model, both $\gamma_k$ and $\beta_k$ ($\forall k\in \mathcal{E}$) should be nonzero since the matrix $A$ is not strictly stable for this case. As a result, we use a small positive value $\xi$ for $\gamma$ instead of 0. By recalling the LMI (\ref{LMI}) with multiple nonlinearities, we can obtain a common quadratic Lyapunov function of dynamic system (\ref{swing}) which is valid for all potential post-fault equilibia with system trajectories inside the polytope $\mathcal{P}$. 
For the case in power grids, the positive definite matrix $P$ in LMI (\ref{LMI}) is of $(|\mathcal{N}|+|\mathcal{N}_G|) \times (|\mathcal{N}|+|\mathcal{N}_G|)$; the matrices $A$, $B$, and $C$ are defined in Subsection IV-A; $\beta=\emph{\emph{diag}}(\beta_1,\dots,\beta_k)$, and $\gamma=\emph{\emph{diag}}(\xi,\dots,\xi)$; the identity matrix $I$ is of $|\mathcal{E}|\times|\mathcal{E}|$. To make the expression of fault-cleared state consistent with the Lur'e reformulation, we have 
\begin{equation} \label{clearing.point}
x(t_c)=[\delta(t_c)-\theta^\prime,\,\omega(t_c)-1]^T.
\end{equation}
where $\delta(t_c)$ and $\omega(t_c)$ are given in (\ref{approximation}).

\subsection{A Novel TSCOPF Framework for Power Systems}

Based on the definitions and  preparations given in Subsections IV.A-B, we develop a novel TSCOPF framework for power grids by applying the stability-constraint optimization framework (\ref{framework}) proposed in Section III. The developed TSCOPF model can be expressed as
\begin{subequations} \label{TSCOPF}
\begin{align}
\min_{p^{G},W_{min}}\  &f =\sum_i (a_{i1}(p_i^G)^2+a_{i2}p_i^G)- \epsilon W^{\min} \\
  \mathrm{s.t.}\quad&(\ref{approximation}),\,(\ref{energy}),\, (\ref{clearing.point}),\, \textrm{and} \nonumber\\
  &\mathbf{\emph{p}}^G-\mathbf{\emph{p}}^L-g^p (\mathbf{\emph{V}},\mathbf{\theta},Y) = 0  \\
  &\mathbf{\emph{q}}^G-\mathbf{\emph{q}}^L-g^q (\mathbf{\emph{V}},\mathbf{\theta},Y) = 0 \\
    &S(\mathbf{\emph{V}},Y) \le \overline{S}\\
   &\underline{V},\underline{E}[\theta] \le V,E[\theta] \le \overline{V},\overline{E}[\theta] \\
  &\mathbf{\emph{p}}^{G\prime}-\mathbf{\emph{p}}^{L\prime}-g^p (\mathbf{\emph{V}}^\prime,\mathbf{\theta}^\prime,Y^\prime) = 0  \\
  &\mathbf{\emph{q}}^{G\prime}-\mathbf{\emph{q}}^{L\prime}-g^q (\mathbf{\emph{V}}^\prime,\mathbf{\theta}^\prime,Y^\prime) = 0 \\
  &S(\mathbf{\emph{V}}^\prime,Y^\prime) \le \overline{S} \\
  &\underline{V},\underline{E}[\theta] \le V^\prime,E[\theta^\prime] \le \overline{V},\overline{E}[\theta] \\
  &\underline{p}^{G},\underline{q}^{G} \le p^{G},q^{G} \le \bar{p}^{G},\bar{q}^{G} \\
  &\begin{cases}
W^{\min} \le \frac{\pi (\pi + 2\theta_{ij}^\prime+ 2\alpha_{ij})}{C_i^T
P^{-1} C_i} \\ 
W^{\min} \le \frac{\pi (\pi - 2\theta_{ij}^\prime - 2\alpha_{ij})}{C_i^T
P^{-1} C_i} 
\end{cases},
\end{align}
\end{subequations}
where constraints (18b)-(18e) and (18f)-(18i) denote the power flows of pre- and post-fault systems respectively. Constraints (18k) is the specialization of its general form $\Psi$ for power systems. Constraints (18f)-(18i) are redundant if the system restores to the original topology after the fault is cleared.

Ignoring the intermediate variable $x(t_c)$ in (\ref{energy}) and (\ref{clearing.point}), all state variables of (\ref{TSCOPF}) are in the steady-state domain. Consequently, it is reasonable to consider the region specified by constraints(\ref{energy}), (\ref{clearing.point}), and (18k) as the projection of the transient stability region with respect to a given fault onto the steady-state domain. There exist some convex relaxations for the power flows (18b)-(18c) and (18f)-(18h) \cite{Low}. Since the $K_i$ in (\ref{approximation}) has the same form as the power flow equations, it is easy to extend the convex relaxations of power flows to constraint (\ref{approximation}). We can obtain a convex TSCOPF model by replacing the nonconvex constraints (\ref{approximation}), (18b)-(18c) and (18f)-(18h) with their convex relaxations. 
 
\textbf{Remark 2.} The TSCOPF framework (\ref{TSCOPF}) is valid for considering the cascading failures/attacks \cite{Kinney} or multiple faults with additional computational burden. In summary, the stability criterion of the proposed optimization framework is that the fault-cleared point given by (\ref{clearing.point}) stays inside the post-fault invariant set described by (\ref{energy}) and (18k). For the cases with cascading failures or attacks, pre- and post- fault equilibria and the invariant set can calculated in exactly the same way as in the single-fault cases. The only difference is that we need to calculate the fault-cleared point using a piecewise Taylor expansion. 

\section{Case Study}
\subsection{Introduction to the Test System}
\begin{figure}[h!]
\centering
\includegraphics[width=0.48\textwidth]{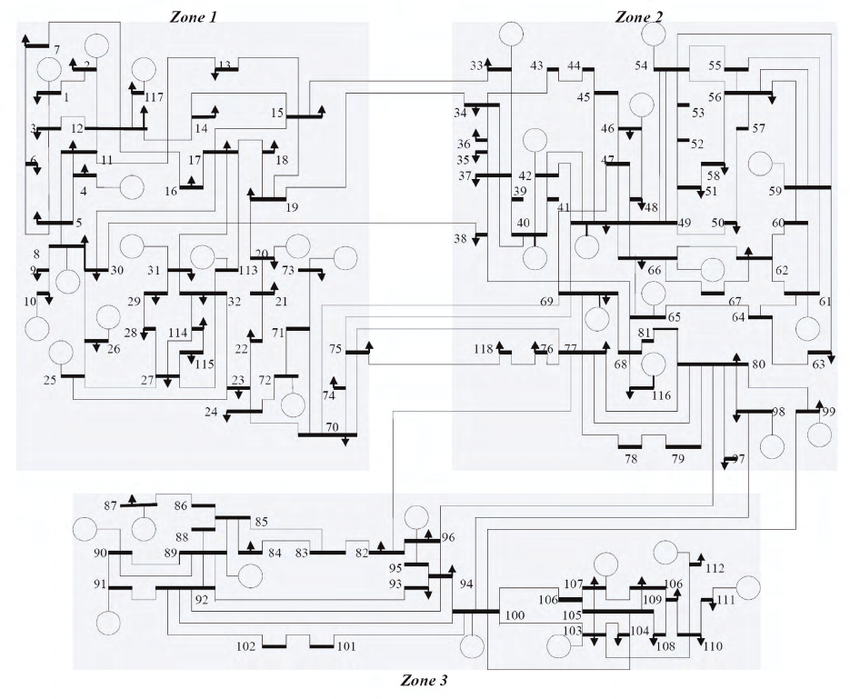}
  \caption{IEEE 118-bus test system.}
  \label{fig:ieee118}
\end{figure}

We test the proposed TSCOPF framework (\ref{TSCOPF}) on one of the most commonly used test transmission grids, the IEEE 118-bus system, which consists of 19 generators, 35 synchronous condensers, 177 lines, 9 transformers, and 91 loads \cite{IEEE118} as shown in Fig. \ref{fig:ieee118}. The dynamic data, i.e. generator moment of inertia $m_i$ ($i \in \mathcal{N}_G$) and the damping coefficient $d_i$ ($i \in \mathcal{N}_G \cup \mathcal{N}_L$), comes from \cite{DynamicData}. The system is stable under a wide range of transient faults with the original load profile. To create two heavy-loaded test cases, we scale up the demand at each load bus by factors of 1.6 and 1.9 respectively. 

\subsection{Case Design and Results}
The results of two scenarios are compared to evaluate the effectiveness of the proposed TSCOPF framework. In the first scenario, we solve the problem (\ref{OPF}) to obtain a solution of the original OPF while problem (\ref{TSCOPF}) is considered in the second scenario. Note that constraints (\ref{clearing.point}) and (18b)-(18e) are related to a specific fault. In this case study, we consider a symmetric bus-to-ground fault at Node 8 since a symmetric fault is the most severe type of fault involving the largest current.
\begin{equation} \label{OPF}
\begin{split}
 \min_{p^{G}}\  &f =\sum_i (a_{i1}(p_i^G)^2+a_{i2}p_i^G) \\
\mathrm{s.t.}\  &(18b)-(18e),  \; \textmd{and} \, (18j)
\end{split}
\end{equation}

\begin{table}[b]
\centering
\caption{Optimal Solutions}
\label{Table1}
\begin{tabular}{cccc}
\hline\hline
Load factor          & Scenario                                                   & \begin{tabular}[c]{@{}c@{}}Objective value\\ (k Dollars)\end{tabular} & \begin{tabular}[c]{@{}c@{}}CPU time\\ (secs)\end{tabular} \\ \hline
\multirow{2}{*}{1.6} & \begin{tabular}[c]{@{}c@{}}Original OPF\end{tabular} & 2341.0998                                                             & 8.986                                                     \\
                     & TSCOPF                                                 & 2350.5974                                                             & 15.815                                                    \\ \hline
\multirow{2}{*}{1.9} & \begin{tabular}[c]{@{}c@{}}Original OPF\end{tabular} & 2895.7271                                                             & 9.571                                                     \\
                     & TSCOPF                                                 & 2909.8168                                                             & 18.298                                                    \\ \hline\hline
\end{tabular}
\end{table}

 The accuracy of the fault-on trajectory approximation (\ref{approximation}) relies on the fault duration. To assess the sensitivity of the proposed TSCOPF model to the fault clearing time, we consider fault clearing times of 0.167s (i.e. 10 cycles) and 0.100s (i.e. 6 cycles) for the 1.6 times and 1.9 times load cases respectively. The optimal solutions, including the required CPU times, of the original OPF and the TSCOPF are tabulated in Table \ref{Table1}. Problems (\ref{OPF}) and (\ref{TSCOPF}) are solved by the nonlinear solver IPOPT (version 3.12.4) \cite{ipopt} through the optimization package JuMP in Julia (version 0.5.2) \cite{JuMP}. Before that, we obtain a uniform quadratic Lyapunov function for (\ref{TSCOPF}) using a MATLAB toolbox YALMIP \cite{yalmip} by calling the SDP solver MOSEK (version
7.1.0.34) \cite{mosek}. A MAC computer with a 64-bit Intel i7 dual core CPU at 2.40 GHz and 8 GB of RAM was used to solve the optimization cases.

\begin{figure}[h]
\centering
\subfloat[Original OPF]{\includegraphics[clip,width=\columnwidth]{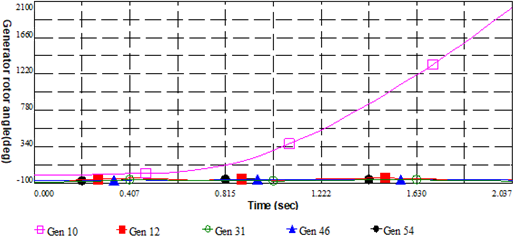}}

\subfloat[TSCOPF]{%
  \includegraphics[clip,width=0.99\columnwidth]{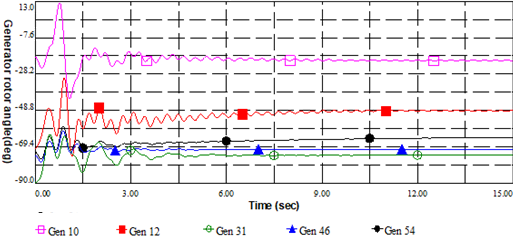}%
}
  \caption{System Responses of the case with load factor of 1.6.}
  \label{fig:1.6load}
\end{figure}

\begin{figure}[h]
\centering
\subfloat[Original OPF]{\includegraphics[clip,width=\columnwidth]{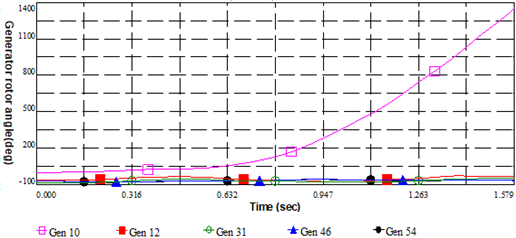}}

\subfloat[TSCOPF]{%
  \includegraphics[clip,width=0.99\columnwidth]{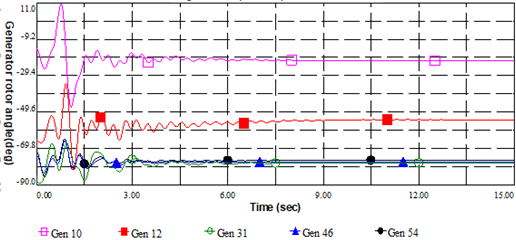}%
}
  \caption{System Responses of the case with load factor of 1.6.}
  \label{fig:1.9load}
\end{figure}

With the initial solutions (i.e. the optimal solutions) obtained in the previous step, transient stability analysis is conducted using DSA Toolbox \cite{dsatool} in a Windows computer with a 64-bit Intel i7-7700 4 cores CPU at 3.60 GHz and 16 GB of RAM. The system responses of the two scenarios are plotted in Fig. \ref{fig:1.6load} and Fig. \ref{fig:1.9load} respectively. 

\subsection{Analysis}
Based on the numerical results obtained in Subsection IV-B, we have the following observations:

(1) The system remains stable after the fault is cleared if it is operated at the equilibrium obtained by solving the proposed TSCOPF. The solution of the original OPF can not guarantee the stability even though the bus-to-ground fault at Bus 8 is cleared very fast.

(2) One can obtain a stability-guaranteed solution by solving the TSCOPF with the cost of a higher (however not significantly higher) computational time, since only a limited number of nonlinear and nonconvex constraints need to be added to the original OPF to construct the TSCOPF. Compared with the existing methods, this is one of the most significant advantages of the proposed TSCOPF framework. To be more precisely, the added constraints include a set of power flow equations (18f)-(18i) with respect to the post-fault equilibrium, the approximation of fault-on trajectories (\ref{clearing.point}) which has a similar form as the power flows with a much smaller size, and some linear or convex constraints (18k) and (\ref{energy}).

(3) The feasible set of TSCOPF is a subset of that of the original OPF. Consequently, it is direct to know that the cost, namely the optimal objective value, of the TSCOPF is higher than that of the original OPF. As shown in Table \ref{Table1}, the costs are increased by less than 5\%, which can be claimed acceptable. Due to the convex relaxation (18k), the conservativeness of the proposed approach is effectively reduced.

\section{Conclusions and Future Work}
This paper proposes a stability-constrained optimization framework for a type of nonlinear systems whose dynamics can be described by a Lur'e system. Unlike the existing methods which are based on either DAE-discretization or iterative algorithms where an independent stability assessment is required for each iteration, the introduced framework is developed based on Lyapunov stability theories. One of the primary advantages of the developed approach is more computational tractable than the existing methods. To illustrate the application values of the proposed framework, it has been successfully applied in power grids to develop a novel TSCOPF model.

The numerical study on the IEEE-118 test system demonstrates that the prosed TSCOPF framework can effectively obtain a stability-guarantee optimal solution with acceptable computational burden. However, There are some directions that can be pursued to push the introduced stability-constrained optimization framework to the online application level. First and foremost, methods of constructing uniform quadratic Lyapunov functions for the nominal system of (\ref{Lure}) with less conservativeness need to be explored. Although the quadratic form Lyapunov fucntions are computationally effective, they may be conservative for stability assessment of many dynamic systems. It is necessary to customize the procedure introduced in Subsection III-B for specific dynamic systems to construct quadratic Lynapunov functions with less conservativeness. Actually, we are currently exploring the possibility of reducing the conservativeness by co-optimizing the sector bound and the coefficient matrix $P$.

The non-convex nature of many Lur'e type systems prevents the application of powerful convex optimization approaches. It is valuable to explore effective convex relaxations to make the stability-constrained optimization framework more computationally tractable for the purpose of online application. From the perspective of power grids, the dynamic system model considering higher-order generator and load models is no longer Lur'e-type and has multi-variate nonlinear terms. It is also valuable to extend the proposed framework to the cases with complex dynamic models.

As mentioned in Remark 2, with additional computational burden, the proposed TSCOPF (\ref{TSCOPF}) is valid for the cases taken into account multiple faults. One of our resent research demonstrates that dynamic response of a microgrid to any fault in the network depends mostly on the type of fault (1-phase, 3-phase etc.) and the fault-clearing time, while having only weak dependence on the fault location and post-fault network topology. This property can be naturally leveraged for developing an effective TSCOPF model for microgrids with multiple faults considered in one scenario.

\appendices

\section{Proof of the Proposition}
 The derivative of Lyapunov function $W$ along the trajectories of the nominal system in (\ref{Lure}) is given by
\begin{align} \label{Wdot}
\dot{W}(x)&=\dot{x}^TPy+y^TP\dot{x} \nonumber \\
&=x^T(A^TP+PA)x+x^TPB\phi+\phi^TB^TPx \nonumber \\
&= \left[\begin{array}{c}
x \\
\phi
\end{array} \right]^T \left[\begin{array}{cc}
A^TP+PA&PB \\
B^TP& 0
\end{array} \right]\left[\begin{array}{c}
x \\
\phi
\end{array} \right] \nonumber . 
\end{align}
Denote the above block matrix as $BM1$. Condition (\ref{lemma}) in the Lemma is equivalent to 
\begin{align} 
\left[\begin{array}{c}
x \\
\phi
\end{array} \right]^T \left[\begin{array}{cc}
\gamma\beta C^TC & -\frac{(\gamma+\beta)}{2}C^T \\
-\frac{(\gamma+\beta)}{2}C& I
\end{array} \right]\left[\begin{array}{c}
x \\
\phi
\end{array} \right] < 0 \nonumber . 
\end{align}
Denote the above block matrix as $BM2$. By observing the LMI in (\ref{LMI}), we realize that $LMI(\ref{LMI})=BM1-BM2 \preceq 0$, which means $BM1 \preceq BM2 \prec 0$ ($\forall y \in \mathcal{P}-\{0\}$). As a result, $\dot{W}(x) < 0$.

\section{Proof of the Lemma}

Suppose $x^\dagger$ is a point on one of the edges of the polytope $\mathcal{P}$. Due to the definition of $W^{\min}$, we have $W(x^\dagger) \geq W^{\min}$. Hence, the system is not able to evolve from the fault clearing point $x(t_c)$ to $x^\dagger$ since $W(x(t_c)) \le W^{\min} \le W(x^\dagger)$. The Lyapunov function value $W(y)$ can only decrease along the system trajectory since $\dot{W} \le 0$ within polytope $\mathcal{P}$. as a result, the system trajectory will stay within the region $\Omega$ or even converge to the origin as $t$ goes to infinity.

\section{Proof of Theorem 1}

Assume that $\hat{s}=[\hat{x}^*$; $\hat{z}$; $\hat{W}^{\min}]$ is an optimal solution of (SL-SCO):\\
\textbf{i.} $\hat{s}$ is feasible to (BL-SCO). 

First, we will show that any optimal solution of (SL-SCO) is optimal to problem (\ref{Wmin}). Suppose that $\hat{W}^{\min}$ does not make equal sign hold in any of (\ref{Concave}), which means $\hat{s}$ is not feasible to (\ref{eq4}). It suffices to show there exists another feasible solution of (SL-SCO), $\bar{\bar{s}}=[\hat{x}^*$; $\hat{z}$; $\hat{W}^{\min}+\Delta W]$, where $\Delta W$ is an arbitrarily small positive value. We have  
\begin{equation}
F(\bar{\bar{s}})-F(\hat{s})=-\epsilon \Delta W \le 0, \nonumber
\end{equation}
which contradicts the optimality of $\hat{s}$. In other words, for any solution $\hat{s}$ in which no equal sign holds in any inequality of (\ref{Concave}), one can always choose a sufficiently small positive value $\Delta W$ to construct another feasible solution $\bar{\bar{s}}$ of (SL-SCO) with a smaller/better objective value until equal sign holds in one of the equalities of (\ref{Wmin}). Such an solution is exactly the optimal solution of problem (\ref{Wmin}). Optimal solution $\hat{s}$ of (SL-SCO) being optimal to (\ref{Wmin}) implies that it is feasible to (BL-SCO) since (BL-SCO) is a bilevel optimization problem with (\ref{Wmin}) as the lower-level problem. \\
\textbf{ii.} $\hat{s}$ is also optimal to (BL-SCO). 

Suppose $\hat{s}$ is not a locally optimal solution of (BL-SCO), then there exist a feasible solution of (BL-SCO), $\bar{\bar{s}}=[\bar{\bar{x}}^*$; $\bar{\bar{z}}$; $\bar{\bar{W}}^{\min}]$, that is in the vicinity of $\hat{s}$ satisfying
\begin{equation}
f(\bar{\bar{x}}^*, \bar{\bar{z}})- f(\hat{x}^*,\hat{z}) \le 0. \nonumber
\end{equation}
It is straightforward to show that all solutions which are feasible to subproblem (\ref{Wmin}) will also satisfy constraint (\ref{Concave}). Thus, $\bar{\bar{s}}$ is also feasible to (SL-SCO) and satisfies
\begin{equation} \label{Contradition}
F(\bar{\bar{s}})-F(\hat{s})=f(\bar{\bar{x}}^*,\bar{\bar{z}})- f(\hat{x}^*,\hat{z}) + \epsilon (\hat{W}^{\min}-\bar{\bar{W}}^{\min}) \le 0.
\end{equation}
Note that $\epsilon$ is an arbitrarily small value. Hence, it is reasonable to assume that the term $\epsilon$($\hat{W}^{\min}$ - $\bar{\bar{W}}^{\min}$) is not comparable to ($f$($\bar{\bar{x}}^*,\bar{\bar{z}}$) - $f$($\hat{x}^*,\hat{z}$)), which means ($F$($\bar{\bar{s}}$) - $F$($\hat{s}$)) has the same sign as ($f$($\bar{\bar{x}}^*,\bar{\bar{z}}$) - $f$($\hat{x}^*,\hat{z}$)). Condition (\ref{Contradition}) contradicts the optimality of $\hat{s}$ to (SL-SCO). Namely, $\hat{s}$ is an local minimum of (BL-SCO) if it is an local minimum of (SL-SCO).

So far, the theorem has been proved by contradiction.

\section{Proof of Theorem 2}

For the sake of convenience, we replace the terms $C_i^Tx^*$ and $C_i^TP^{-1}C_i$ with $X$ and $1/\lambda$ respectively. The notation $CONV(A)$ means the convex hull of set $A$.\\
\textbf{i.} $CONV(\psi) \subseteq \Psi$.

For any $X \in [\underline{X},0]$, we have 
\begin{align}
\psi &=\{(X,W^{\min})|0 \le W^{\min} \le \lambda(X+\Delta l)^2\} \nonumber \\
\Psi &=\{(X,W^{\min})|0 \le W^{\min} \le \lambda((\underline{X}+2\Delta l)X+\Delta l^2)\}. \nonumber
\end{align}
It is direct to know that $\psi \subseteq \Psi$ since $(X+\Delta l)^2 \le (\underline{X}+2\Delta l)X+\Delta l^2)$. Similarly, for any $X \in [0,\overline{X}]$, we have the same conclusion that $\psi \subseteq \Psi$, which means $\Psi$ is a convex relaxation of $\psi$. Since convex hull is defined as the intersection of all convex relaxations of a non-convex set \cite{Li2}, we have $CONV(\psi) \subseteq \Psi$.\\
\textbf{ii.} $CONV(\psi) \supseteq \Psi$.

If a linear inequality is valid for a given set $\Omega_A$, it will also be valid for any subset of $\Omega_A$. Note that “a linear inequality is valid for a set” means the inequality is satisfied by all its feasible solutions \cite{Validcut}. On the other way round, according to the properties of supporting hyperplanes \cite{Boyd2}, $\Omega_B$ is said to be a subset of $\Omega_A$ if $\Omega_A$ is convex and any valid linear inequality of $\Omega_A$ is also valid for $\Omega_B$ \cite{Li3}. Let $s=[X$; $W^{\min}]$ and suppose that $\alpha s \geq \beta$ is any valid linear cut for $CONV(\psi)$, this cut should be also valid for all the points in $\psi$. To prove that $CONV(\psi) \supseteq \Psi$ (i.e. $\Psi$ is a subset of $CONV(\psi)$), we try to show that $\alpha s \geq \beta$ is valid for all the edges of $\Psi$.

The convex set $\Psi$ has five edges of which the formulations are given as
   \[   
    \Psi_1 = \left\{(X,W^{\min}) \left| \begin{array}{lr}
W^{\min}= \lambda ((\overline{X} - 2\Delta l)X+\Delta l^2) \\ 
W^{\min} \le \lambda ((\underline{X} + 2\Delta l)X+\Delta l^2) \\
\underline{X} \le X \le \overline{X}
\end{array}\right. \right\}, 
  \]
     \[   
    \Psi_2 = \left\{(X,W^{\min}) \left| \begin{array}{lr}
W^{\min} \le \lambda ((\overline{X} - 2\Delta l)X+\Delta l^2) \\ 
W^{\min} = \lambda ((\underline{X} + 2\Delta l)X+\Delta l^2) \\
\underline{X} \le X \le \overline{X}
\end{array}\right. \right\}, 
  \]
     \[   
    \Psi_3 = \left\{(X,W^{\min}) \left| \begin{array}{lr}
W^{\min} \le \lambda ((\overline{X} - 2\Delta l)X+\Delta l^2) \\ 
W^{\min} \le \lambda ((\underline{X} + 2\Delta l)X+\Delta l^2) \\
 X = \underline{X}
\end{array}\right. \right\}, 
  \]
     \[   
    \Psi_4 = \left\{(X,W^{\min}) \left| \begin{array}{lr}
W^{\min} \le \lambda ((\overline{X} - 2\Delta l)X+\Delta l^2) \\ 
W^{\min} \le \lambda ((\underline{X} + 2\Delta l)X+\Delta l^2) \\
X = \overline{X}
\end{array}\right. \right\}, 
  \]
     \[   
    \Psi_5 = \left\{(X,W^{\min}) \left| \begin{array}{lr}
W^{\min} = 0 \\
\underline{X} \le X \le \overline{X}
\end{array}\right. \right\}. 
  \]
  
As an example, we show that the cut $\alpha s \geq \beta$ is valid for edge $\Psi_1$ in this paragraph. It is easy to verify that the two points $s_1=(0,\lambda \Delta l^2)$ and $s_2=(\overline{X},\lambda(\overline{X}-\Delta l)^2)$ are located in both $\psi$ and $\Psi_1$. That means the cut $\alpha s \geq \beta$ is valid for these two points and we have $\alpha s_1 \geq \beta$ and $\alpha s_2 \geq \beta$. Let $\hat{s}=(\hat{X},\hat{W}^{min})$ denote any given point in set $\Psi_1$. For any given $\hat{s}$, there exists a value $c$ $(0 \le c \le 1)$ satisfying $\hat{s}=cs_1+(1-c)s_2)$. It suffices to verify this statement by substituting $\hat{s}=((1-c)\overline{X},c\lambda \Delta l^2 +(1-c)\lambda(\overline{X}-\Delta l)^2)$ into the first equation in $\Psi_1$. As a result, we have
\begin{equation}
\alpha \hat{s}= c\alpha s_1+(1-c) \alpha s_2 \geq c\beta +(1-c)\beta = \beta, \nonumber
\end{equation}
which means the linear cut is also valid for any given point in $\Psi_1$ and, consequently, valid for $\Psi_1$.

Using the same method, the readers are able to prove that the linear cut $\alpha s \geq \beta$ is valid for all the other four edges of $\Psi$ and, consequently, valid for the whole convex set $\Psi$. Hence, $CONV(\psi) \supseteq \Psi$.\\
\textbf{iii.} $\Psi = CONV(\psi)$.

$CONV(\psi) \subseteq \Psi$ and $CONV(\psi) \supseteq \Psi$ together imply $\Psi = CONV(\psi)$.

\ifCLASSOPTIONcaptionsoff
  \newpage
\fi

\end{document}